\documentclass[12pt]{article}
\usepackage[tbtags]{amsmath}
\usepackage{epsfig,amstext,amssymb,amsthm,latexsym}
\catcode`\@=11
\@addtoreset{equation}{section}
\catcode`\@=12

\usepackage{amssymb,color}
\definecolor{c20}{rgb}{0.,0.7,0.}
\definecolor{c30}{rgb}{0.,0.,1.}
\definecolor{c40}{rgb}{1,0.1,0.7}
\definecolor{c50}{rgb}{1,0,0}

\allowdisplaybreaks[4]
\newcommand{\nwc}{\newcommand}
\nwc{\COM}[1]{}

\nwc{\vs}[1]{\vskip #1 cm}

\newtheorem{theo}{Theorem}[section]
\newtheorem{sat}[theo]{Proposition}
\newtheorem{de}[theo]{Definition}
\newtheorem{lem}[theo]{Lemma}

\newtheorem{korr}[theo]{Corollary}

\newtheorem{exxa}[theo]{Example}

\newcommand{\nelem}[1]{{Lemma \ref{#1}}}

\newcommand{\netheo}[1]{{Theorem \ref{#1}}}

\newcommand{\kb}[1]{\boldsymbol{#1}}
\newcommand{\vk}[1]{\kb{#1}}

\def\FRE{\mbox{Fr\'{e}chet }}
\def\a{\vk{a}}
\def\a{\alpha}

\newcommand{\ve}{\varepsilon}

\newcommand{\abs}[1]{\lvert #1 \rvert}
\newcommand{\Abs}[1]{ \Bigl \lvert #1 \Bigr \rvert}

\newcommand{\E}[1]{\mbox{\rm$\vk{E}$}\{#1\}}
\newcommand{\eg}[1]{\mbox{\rm$\vk{E}$}\biggl\{#1 \biggr\}}

\newcommand{\pk}[1]{\mbox{\rm$\vk{P}$} \{#1\} }

\newcommand{\pb}[1]{\mbox{\rm$\vk{P}$}\Bigl \{#1 \Bigr \}}

\newcommand{\R}{\!I\!\!R}

\newcommand{\inr}{\in \R}

\newcommand{\limit}[1]{\lim_{#1 \to   \infty}}

\newcommand{\equaldis}{\stackrel{d}{=}}

\newcommand{\BQN}{\begin{eqnarray}}
\newcommand{\EQN}{\end{eqnarray}}
\newcommand{\BQNY}{\begin{eqnarray*}}
\newcommand{\EQNY}{\end{eqnarray*}}
\newcommand{\BS}{\begin{sat}}
\newcommand{\ES}{\end{sat}}
\newcommand{\BL}{\begin{lem}}
\newcommand{\EL}{\end{lem}}
\newcommand{\BT}{\begin{theo}}
\newcommand{\ET}{\end{theo}}
\newcommand{\BK}{\begin{korr}}
\newcommand{\EK}{\end{korr}}
\newcommand{\BD}{\begin{de}}
\newcommand{\ED}{\end{de}}
\newcommand{\BIT}{\begin{itemize}}
\newcommand{\EIT}{\end{itemize}}
\newcommand{\BDI}{\begin{description}}
\newcommand{\EDI}{\end{description}}

\newcommand{\QED}{\hfill $\Box$}

\newcommand{\IF}{\infty}
\def\kal#1{{\cal{ #1}}}
\def\fracl#1#2{\biggr( \frac{#1}{#2} \biggl) }

\newcommand{\prooftheo}[1]{ \textsc{Proof of Theorem} \ref{#1} }

\newcommand{\prooflem}[1]{\textsc{Proof of Lemma} \ref{#1}}

\newcommand{\BEX}{\begin{exxa}}
\newcommand{\EEX}{\end{exxa}}

\def\st{\stackrel}

\def\omgF{{r_F}}
\def\toomgF{ \uparrow \omgF}

\begin{document}

\centerline{\Large Asymptotics of Random Contractions}

        \vskip 1.8 cm
        \centerline{\large ENKELEJD HASHORVA \footnote{
       Department of Actuarial Science, Faculty of Business and Economics,
University of Lausanne, B\^{a}timent Extranef, UNIL-Dorigny, 1015 Lausanne, Switzerland,
E-mail: Enkelejd.Hashorva@unil.ch}, ANTHONY G. PAKES \footnote{School
of Mathematics and Statistics, University of Western Australia, 35
Stirling Highway, Crawley, W. A., 6009, Australia, E-mail:
pakes@maths.uwa.edu.au},
 and QIHE TANG\footnote{
Department of Statistics and Actuarial Science, The University of
Iowa, 241 Schaeffer Hall, Iowa City, IA 52242-1409, United States,
E-mail: qtang@stat.uiowa.edu} }

        \vskip 1.4 cm

{\bf Abstract:} In this paper we discuss the asymptotic behaviour of
random contractions $X=RS$, where $R$, with distribution function
$F$, is a positive random variable independent of $S\in (0,1)$.
Random contractions appear naturally in insurance and finance. Our
principal contribution is the derivation of the tail asymptotics of
$X$ assuming that  $F$ is in the max-domain of attraction of  an
extreme value distribution and the distribution function of $S$
satisfies a regular variation property. We apply our result to
derive the asymptotics of the probability of ruin for a particular
discrete-time risk model. Further we quantify in our asymptotic
setting the effect of the random scaling on the Conditional Tail
Expectations, risk aggregation, and derive the joint asymptotic
distribution of linear combinations of random contractions.



\section{Introduction}
Let $R,S$ be two independent random variables with $R>0, S\in (0,1)$
almost surely, and define $X= RS$. The random variable $X$ is a
random contraction of $R$ via $S$. Random contractions or random
scalings are common in insurance and finance applications. Typically
$R$ models a random payment whereas $S$ is a random discount factor.
Several authors have studied random contractions in quite different
contexts. Some recent contributions dealing with distributional and
asymptotic properties of random contractions are Kotz and Nadarajah
(2000), Galambos and Simonelli (2004),  Gomes et al. (2004), Maulik
and Resnick (2004), Tang  and Tsitsiashvili  (2003, 2004), Jessen
and Mikosch (2006), D'Auria  and Resnick (2006, 2008), Tang (2006,
2008), Denisov  and Zwart (2007), Pakes and Navarro
(2007), Resnick (2007), Beutner and Kamps (2008), Charpentier and
Segers (2007, 2009), Hashorva (2008, 2009, 2010), Hashorva and Pakes
(2010), Liu and Tang (2010).

Our main goal in this paper is to investigate the tail asymptotics
of random contractions when the tail asymptotics of $R$ is known. An
important motivation for this investigation is the fact that in
insurance and finance applications assumptions often are made on the
tail behaviour of a random payment modeled by $R$. If $S$ represents
the random discount factor applicable to the interval from the
present to the payment time, then $RS$ is the present value of the
later payment  $R$. In cases where the distribution function $G$ of
$S$ is  unknown, it is of some interest to know how the tail
behaviour of the random contraction $X=RS$ is determined by
the corresponding asymptotic behaviours of the factors. One possible application is
to approximating the Value at Risk in the presence of discounting,
given  information about the Value at Risk before discounting.

Without going into mathematical details, we mention briefly the main contributions of this paper: \\
a) Under the assumption that the distribution  function $F$ of $R$ is in the max-domain of attraction of some
univariate extreme value distribution we obtain the asymptotic behaviour of $\pk{X > u}$ as  $u$ tends to the upper endpoint of
$F$, provided that $\overline G=1-G$ satisfies a regular variation property (\netheo{T1} below). \\
b) We determine corresponding results for the density function of
$X$ assuming a regular variation property for  the density
function of $S$, and additional regularity properties in some cases.\\
c) We present four applications: c1) First we derive the asymptotics
of the ruin probability for a particular discrete-time ruin problem,
c2) then we discuss briefly the asymptotics of Conditional Tail
Expectation in the random contraction framework, c3) and we obtain
asymptotic expansions for the aggregation of two contractions, which
lead to novel asymptotic characterisation of bivariate elliptical
distributions, c4) finally we show the asymptotic independence of
certain bivariate random contractions.

As mentioned above, we assume that a generic scaling factor $S$ with
distribution function $G$ takes values in $(0,1)$. In addition, we
assume that $\overline G(1-y)$ is regularly varying at zero, the
above mentioned regular variation property. However, the reader will
easily appreciate that the scaling factors can be multiplied by a
positive constant and hence can function as inflation or deflation
factors. Since our results can easily be adjusted for this
contingency, we will say no more about it beyond the closure
property in \nelem{000}.

It is interesting that under the setup of this paper the asymptotic
tail behaviours of $R$ and $X$ are very similar. In particular,
membership of a max-domain of attraction is insensitive to the
distribution of bounded discount factors.

Our main results are presented in Section 3 followed by the applications in Sections 4.
The proofs of all the results are relegated to Section 5.

\section{Maximal Domains of Attraction}
In this short section we present some details on max-domains of
attraction. The distribution function $F$ belongs to the max-domain of
attraction of a univariate extreme value distribution function $N$,
written $F\in {\rm MDA}(N)$,  if
 \BQN \label{eq:LL}
 \limit{n}\sup_{x\inr} \Abs{F^n(a_nx+ b_n)- N(x)}= 0
 \EQN
holds for some constants $a_n>0,b_n\inr,n\ge 1$. See e.g., Reiss
(1989), Embrechts et al.\ (1997), Falk et al.\ (2004), De Haan and
Ferreira (2006), or Resnick (2008) for more details on univariate
max-domains of attraction. Only three choices for $N$ are possible,
namely the \FRE distribution, the Gumbel distribution, or the
Weibull distribution. We denote the corresponding distribution
functions by $\Phi_\gamma$,  $\Lambda$, and $\Psi_\gamma$,
respectively, where $\gamma>0$ indexes members of the \FRE and
Weibull families.

The functional form of the \FRE distribution function is
$\Phi_\gamma(x)=\exp(- x^{-\gamma})$, $x>0$. If  $F\in {\rm
MDA}(\Phi_\gamma)$,  then \eqref{eq:LL} with
$N=\Phi_\gamma$ is equivalent to
 \BQN\label{eq:PhiF}
 \limit{u} \frac{\overline{F}(xu)}{\overline{F}(u)}= x^{-\gamma}, \qquad \forall
 x>0.
 \EQN
This means that the survival function $\overline F=1- F$  is regularly varying at infinity with
index $-\gamma$ and further it has an infinite upper endpoint (denoted in the sequel by $r_F$).

The functional form of the standard Gumbel distribution function is
$\Lambda(x)=\exp(-\exp(-x)),x\inr$, and  \eqref{eq:LL} with
$N=\Lambda$  is equivalent to
 \BQN \label{eq:rdfd} \lim_{u \uparrow
\omgF} \frac{\overline{F}(u+x/w(u))}{\overline{F}(u)} =
\exp(-x),\qquad \forall x\inr,
 \EQN
where $w$ is a positive scaling function satisfying
\BQN\label{eq:uv}
 \lim_{u\uparrow \omgF} u w(u)=\infty, \quad \text{and  }
\lim_{u\uparrow \omgF} w(u)(\omgF - u) =\infty \quad \text{if  }
r_F< \infty.
 \EQN
Recall that the scaling function $w$ can be defined asymptotically
via the mean excess function (see e.g., Embrechts et al.\ (1997) or
Resnick (2008)) by
 \BQN \label{eq:reco1} w(u)\sim \frac{1}{\E{R-u
\lvert R>u}}, \qquad u\uparrow r_F.
 \EQN
Throughout this paper he relation $a(u)\sim b(u)$ means that the quotient of
both sides tends to $1$ according to the indicated limit procedure.\\
The functional form of the Weibull distribution function is
$\Psi_\gamma(x)= \exp(-\abs{x}^\gamma)$, $x\le 0$. If $F\in {\rm
MDA}(\Psi_\gamma)$, then $r_F$ is finite and \eqref{eq:LL}  is
equivalent to
 \BQN\label{eq:PsiF}
 \limit{u} \frac{\overline{F}(r_F- x/u)}{\overline{F}(r_F-1/u)} =x^{\gamma}, \qquad \forall x>0.
 \EQN

In some applications it is necessary to admit scaling factors $S\in
(0,c)$, where $c$ is a positive constant. We always assume that
$c=1$.  If $F \in {\rm MDA}(\Psi_\gamma)$ we assume too that
$r_F=1$.  The next lemma explains why these
conventions are not restrictive.\\
\BL \label{000} Let $W$ be a random variable whose
distribution function $F$ satisfies \eqref{eq:LL}. If $ c,p\in
(0,\IF)$, then $cW^p$ has a distribution function in the same
max-domain of attraction as $W$. \EL

We will need the following facts about subexponential distribution
functions. See Embrechts et al.\ (1997, Appendix A3) for
distribution functions supported in $[0,\infty)$, and Borovkov and
Borovkov (2008, p. 13) for the general case.  Let $F^{2*}$ denote
the convolution square of $F$, i.e., the distribution function of
$R+R^*$, where $R^*$ is an independent copy of $R$. Assume that
$r_F=\infty$ in what follows.

 In the case that $F(0-)=0$, we say that $F$ is subexponential, written $F\in {\cal S}_+$, if
 \BQN \label{sep}
 \lim_{u\to\infty} {\overline{F^{2*}}(u) \over \overline F(u)}=2.
 \EQN
 In the case that $F$ is two-sided, i.e., $F(0-)>0$, define $F_+(u)=0$ if $u<0$ and $F_+(u)=F(u)$ if $u\geq0$.
 We say that $F$ is subexponential if $F_+ \in {\cal S}_+$, and we write $F\in {\cal S}$; see Borovkov and Borovkov (2008,
 p.
 14).
 These authors (p. 19) show that it is possible that a two-sided distribution function $F\not\in {\cal S}$ satisfies
 \eqref{sep}, but imposing a condition described below will ensure that $F\in {\cal S}$ is equivalent to \eqref{sep}.

 Say that $F$ belongs to the class of long-tailed distribution functions, written  $F\in {\cal L}$, if
 $$\lim_{u\to\infty} {\overline F(u+y) \over \overline F(u)}=1$$
 for all real $y$. The convergence here is locally uniform with respect to $y$. If $F(0-)=0$, \eqref{sep} implies that $F\in
 {\cal L}$.
 In the two-sided case, if $F\in {\cal L}$, then $\overline{F^{2*}}(u)\sim \overline{F_+^{2*}}(u)$ as $u\to\infty$
 (Borovkov and Borovkov (2008,  Theorem 1.2.4(vi)), and hence  $F\in {\cal S}$ is equivalent to \eqref{sep}.
 These concepts relate to attraction to the Gumbel distribution as follows.

 Still assuming that $r_F=\infty$, assume too that $F\in {\rm MDA}(\Lambda,w)$ with $\lim_{u\to\infty}w(u)=0$. For real
 $y$
 we
 can choose $u$ so large that $|y| \le | x|/w(u)$ and hence conclude from \eqref{eq:rdfd} that $F\in {\cal L}$. Further
 conditions can be given to ensure that $F\in {\cal L}\cap{\cal S}$. This holds if there is a positive constant $\lambda$ such
 that (see Mitra and Resnick (2008, Corollary 2.9))
 \BQN \label{eq:sub:R}
\limit{u} \frac{ [\overline F( \lambda/w(u))]^2}{\overline F(u)}= 0.
 \EQN
Note that if $F(0-)=0$, in view of Corollary 2.5 of Goldie and
Resnick (1988), when $w$ is eventually non-increasing such that
 \BQN
\label{grc}
 \lim_{u\to\infty} {w(u) \over w(tu)}>1
 \EQN
holds for some constant $t>1$, then $F\in \kal{S}_+$.

We derive by our next result a  self-contained proof of the Mitra-Resnick criterion \eqref{eq:sub:R}. \\
\BL \label{interesting}
Let $F\in \mathrm{MDA}(\Lambda ,w)$ with
$r_{F}=\infty $
and $\lim_{u\rightarrow \infty }w(u)=0$ (hence $F\in \mathcal{L}$). Then $%
F\in \mathcal{S}$ if and only if
\BQN \label{Tony}
\lim_{u\rightarrow \infty }{\frac{1}{\overline{F}(u)}}\int_{\lambda
/w(u)}^{u-\lambda /w(u)}\overline{F}(u-y)dF(y)=0
 \EQN
holds for some $\lambda >0$. \EL

By \nelem{interesting}, the Mitra-Resnick criterion follows
immediately, because the integral in \eqref{Tony} is bounded above
by
$$
\overline{F}(\lambda /w(u))\left[ \overline{F}(\lambda /w(u))-\overline{F}%
(u-\lambda /w(u))\right] \leq \left[ \overline{F}(\lambda
/w(u))\right] ^{2}.
$$
Furthermore, \nelem{interesting} implies a generalization of
Goldie's sufficient condition asserting that if $F\in
\mathcal{L}$ has a dominatedly varying right-hand tail, i.e., $%
\limsup_{u\rightarrow \infty
}\overline{F}(u/2)/\overline{F}(u)<\infty $, then $F\in \mathcal{S}$
(see Embrechts et al. (1997, pp. 49, 52)). More precisely, suppose
that $F\in \mathrm{MDA}(\Lambda ,w)$ with $r_{F}=\infty $
and $\lim_{u\rightarrow \infty }w(u)=0$. Then by \nelem{interesting}, $%
F\in \mathcal{S}$ if
\BQN
\lim_{u\rightarrow \infty }{\frac{\overline{F}(u/2)}{\overline{F}(u)}}%
\overline{F}(\lambda /w(u))&=&0
\EQN
for some $\lambda \in (0,\infty )$.

\section{Principal Results}
Let $R$ be a positive random variable with distribution function
$F$, and let $S_1,\dots,S_n$ be mutually independent, and
independent of $R$. Denote by $G_i$ the distribution function of
$S_i$, and assume it is supported on $(0,1)$.  In insurance and
financial contexts $S_i$ represents the random discount factor over
the interval $[i-1,i)$. Then
$$X_n= R \prod_{i=1}^n S_i$$
is the present value of the payment $R$ received at time $n$. Let
$H_n$ denote the distribution function of the random product $X_n$.

If $F\in {\rm MDA}(\Phi_\gamma)$ then, with no further conditions,
it follows from Breiman's lemma (Breiman (1965)) that
\BQN\label{eq:breiman}
 \limit{u} \frac{\pk{X_n > u}}{\pk{R>u}}= \eg{\prod_{i=1}^n S_i^\gamma}
 \EQN
and, in particular, that  $H_n\in {\rm MDA}(\Phi_\gamma)$. See
Jessen and Mikosch (2006), Denisov and Zwart (2007), and Resnick
(2007) for details on Breiman's lemma and for  some of its
generalisations. So we need to consider only the cases where  $F$ is
in the max-domain of attraction of the Gumbel or of the Weibull
distribution. In these cases it turns out that the tail asymptotic behaviour of
each $S_i$ is crucial. Our working assumption for the scaling random
variables is that
 \BQN\label{eq:alp}
 \limit{u} \frac{ \overline G_i(1- x/u)}{\overline G_i(1- 1/u)}= x^{\alpha_i}, \qquad \forall
 x>0
 \EQN
for some $\alpha_i  \in [0,\IF)$. So if $\alpha_i>0$, then $G_i\in
{\rm MDA}(\Psi_{\alpha_i})$.  If $G_i$ possesses a positive density
function $g_i$, and if $\alpha_i\in (0,\IF)$, then  we will assume
the von Mises condition
 \BQN\label{eq:alp:den}
 \limit{u} \frac{ g_i(1- x/u)}{g_i (1- 1/u)} = x^{\alpha_i-1}, \qquad \forall x>0.
\EQN

The following is our first result, in which we denote by $\Gamma(\cdot)$ the Euler gamma function.\\
\BT \label{T1} For $i\le n$  let $S_i$ (with distribution function
$G_i$) be mutually independent scaling random variables, and
independent of $R$ (with distribution function  $F$). Assume that
$F(0-)=0$ and $r_F\in (0,\IF]$, and also
that \eqref{eq:alp} holds for every $i\le n$ with $\alpha_i\in [0,\IF)$. \\
a) If $F\in {\rm MDA}(\Lambda,w)$ with $r_F\in (0,\IF]$, then
 \BQN\label{eq:lem1:1} \pk{X_n > u}\sim \prod_{i=1}^n
\Bigl[\Gamma(\alpha_i+1) \overline G_i\left(1- \frac{1}{u
w(u)}\right)\Bigr] \overline F (u), \qquad u\toomgF.
 \EQN
a1) If, in addition, $r_F=\infty$, $\lim_{u\to\infty} w(u)=0$, and \eqref{Tony} holds, then $H_n \in {\cal S_+}$. \\
b) If $F\in {\rm MDA}(\Psi_\gamma)$ with $\gamma\in [0,\IF)$ (hence
$r_F=1$ by our convention), then
 \BQN\label{eq:lem1:2} \pk{X_n > u}
\sim {\Gamma(\gamma+1) \over \Gamma(\gamma+\sum_{i=1}^n\alpha_i+1)}
\prod_{i=1}^n \left[\Gamma(\alpha_i+1) \overline G_i(u)
\right]\overline F(u) , \qquad u\uparrow  1.
 \EQN \ET

In the light of \netheo{T1}, if $F$ is in the Gumbel or the Weibull
max-domain of attraction and each survival function $\overline G_i(1-u)$ is regularly
varying at $0$, then $X_n$ has a distribution function in the Gumbel
or the Weibull max-domain of attraction, respectively. See also the
proof of \netheo{T1} below. Part a1) of \netheo{T1} shows that the
subexponentiality is preserved under random scaling. Recent results
on the subexponentiality of products are obtained in Tang (2006,
2008), Liu and Tang (2010).

The random contraction $X_n$ possesses a density function $h_n$ if
one of the scaling random variables $S_i$ has a density function
$g_i$, or if $R$ has a density function $f$. In the following
theorem we derive asymptotic approximations of the density function
$h_n$. Part a) is a density version of Breiman's lemma. The
assumption under a1) below places conditions on one of the density
functions $g_i$ but not on $F$ (beyond \eqref{eq:PhiF}). In
particular, $g_i$ can behave like a beta density near the origin,
but it is bounded near unity. It seems that  relaxing this condition
requires restricting the slowly varying factor of $\overline F$. One
possibility is to assume that this factor is normalized slowly
varying, i.e., that the density function $f$ exists and is regularly
varying. We do this in a2), and then we need no conditions on $G_i$.
In parts b) and c) we assume for some $i$ that $g_i$ satisfies
\eqref{eq:alp:den}
with an additional technical condition for b).\\
\BT \label{T2} Let $S_i, i\le n, n\ge 1$, and $R$ (with distribution function $F$) be as in \netheo{T1}. \\
a) If  $F$ satisfies \eqref{eq:PhiF} with $\gamma\in [0,\IF)$, then
\BQN\label{eq:T2:2} \limit{u} \frac{uh_n(u)}{\pk{X_n> u}}=\gamma
\EQN
holds under either of the following conditions:\\
a1)  For some $i\le n$, the scaling factor $S_i$ has a density function such that $yg_i(y)$ is bounded in $(0,1)$, or\\
a2) $F$ has a density function $f$ which is regularly varying at
infinity with index $-(\gamma+1)$ for some $\gamma\ge 0$,
and $\int_0^1 y^{-\epsilon} dG_i(y)<\infty$ for some $\epsilon>0$ and all $i\le n$ provided that $\gamma=0$. \\
In addition, $\lim_{u\to\infty} h_n(u) / f(u)$ exists and equals the limit in \eqref{eq:breiman}.\\
b) Assume that \eqref{eq:alp} holds for all $i\le n$ with
$\alpha_i\in (0,\IF)$, and $F\in {\rm MDA}(\Lambda,w)$ with $r_F\in
(0,\IF]$. Suppose too for some $i\le n$ that $S_i$ has a density
function $g_i$ satisfying \eqref{eq:alp:den} and that there exist
$c\in(0,1)$ and $p_i>0$ such that
 \BQN \label{gbb} \sup_{0<y\le
c}y^{p_i}g_i(y)<\infty.
 \EQN
Then
 \BQN\label{eq:T2:1} \lim_{u\toomgF} \frac{h_n(u)}{w(u) \pk{X_n>
u}}= 1.
 \EQN
c) Suppose that the distribution function $F$ satisfies
\eqref{eq:PsiF} with $\gamma\in [0,\IF)$ (hence, $r_F=1$ by our
convention) and that \eqref{eq:alp} holds for all $i\le n$ with
$\alpha_i\ge 0$. If in addition \eqref{eq:alp:den} holds
for some $i\le n$ with $\alpha_i>0$, then
\BQN\label{eq:T2:2}
\lim_{u \downarrow 0} \frac{ u h_n(1- u)}{\pk{X_n > 1-
u}}=\gamma+\sum_{i=1}^n\alpha_i. \EQN \ET

Under the assumptions of part b) of \netheo{T2} we obtain further
 \BQN \label{eq:vm:h} \lim_{u\toomgF} \frac{h_n(u+x/w(u))}{h_n(u) }=
\exp(-x), \qquad \forall x\inr.
 \EQN
If $F\in {\rm MDA}(\Lambda,w)$ and it  has a density function $f$,
then conditions exist under which
 \BQN\label{eq:vm:f}
 \lim_{u\toomgF} \frac{f(u+x/w(u))}{f(u) }= \exp(-x), \qquad \forall x\inr;
 \EQN
see Resnick (2008). If \eqref{eq:vm:f} holds, then we can derive
\eqref{eq:vm:h} under milder conditions than in \netheo{T2}(b). In
essence, we assume below that all the $g_i$ exist and satisfy a
condition
similar to intermediate regular variation. \\
\BT \label{TT} Let $S_i, i\le n, n\ge 1$, be mutually independent
scaling random variables, with positive density functions $g_i$, and
independent of the random variable $R\geq 0$ which has a
distribution function $F$ with an upper endpoint $r_F\in (0,\IF]$.
Suppose further that $F$ possesses a positive density function $f$
such that \eqref{eq:vm:f} holds. If for  $i\le n$, $s \in (u,r_F)$,
and any measurable function $a:\R^2 \to [0,\IF)$ such that $
\lim_{u\toomgF} a(u,s)=1$ we have \BQN \label{eq:gu}
 \lim_{u\toomgF} \frac{g_i((u/s) a(u,s))}{ g_i(u/s)}=1,
\EQN
then \eqref{eq:vm:h} is satisfied.
 \ET
Note in passing that if $r_F$ is finite and $g_i(1-u)$ is regularly
varying at $0$ with index $\alpha_i -1 \in (0,\IF)$,  then
\eqref{eq:gu} is satisfied.

We present next three   illustrating examples.

{\bf Example 1.} Let $R$ be a random variable with distribution
function $F\in {\rm MDA}(\Lambda, w)$ and upper endpoint $r_F\in
(0,\IF]$, and let
 $S$ be a random variable with beta distribution with positive parameters $\alpha,\beta$.  Since
 \BQN\label{cited later}
\pk{S>1-  u}\sim
\frac{\Gamma(\alpha+\beta)}{\Gamma(\alpha+1)\Gamma(\beta)} u^\alpha,
\qquad u \downarrow 0,
 \EQN
it follows from  \netheo{T1} that the distribution function $H$ of
$RS$ satisfies
 \BQNY \overline H(u)\sim \frac{\Gamma(\alpha+
\beta)}{\Gamma(\beta)} (u w(u))^{-\alpha} \overline F(u), \qquad
u\uparrow r_F.
 \EQNY
If $g$ is the positive density function of $S$, then condition
\eqref{gbb} holds, whence  \eqref{eq:T2:1} implies that the density
function $h$ of $H$ satisfies
 \BQNY h(u)\sim  w(u) \overline{ H}(u)
, \qquad u\uparrow r_F.
 \EQNY
Note in passing that condition \eqref{eq:gu} can be easily checked.

{\bf Example 2.} Under the setup of the previous example suppose
that
 \BQN\label{eq:kotz:Fk}
  \pk{R> u} \sim K u^{q}\exp(-r u^\gamma
), \qquad K>0,r>0,\gamma>0, q\inr
 \EQN
holds for $u\to \IF$.  It follows easily that $F\in {\rm
MDA}(\Lambda,w)$ with $w(u)= r \gamma u ^{\gamma -1},u>0$. Hence
relation \eqref{cited later} implies
 \BQNY \overline H(u) \sim K (r
\gamma  )^{-\alpha} \frac{\Gamma(\alpha+ \beta)}{\Gamma(\beta)}
u^{q- \alpha \gamma }\exp(-r u^\gamma ), \qquad u\to \IF.
 \EQNY
Further, by \eqref{eq:T2:1}, the density function $h$ of $RS$
satisfies
 \BQNY h(u)\sim r \gamma  u ^{\gamma -1} \pk{ RS > u} ,
\qquad u\to \IF.
 \EQNY
Note in passing that condition \eqref{eq:sub:R} does not hold for
any $\lambda \in (0,\IF)$. However, if $\gamma \in (0,1)$, then $R$
and $RS$ have subexponential distributions.

{\bf Example 3.} Let $R$ be a positive random variable with distribution function $F$ satisfying
$$ \overline{F}(u)\sim c_1 \exp(- c_2/(1- u)) , \qquad u \uparrow 1,$$
with $c_1,c_2$ two positive constants. Since for $w(u)= c_2/(1-
u)^2$, $u\in (0,1)$, and any $s\inr$,
$$ \frac{ \overline{F}(u+ s/w(u))}{ \overline{F}(u)}\sim \exp( - c_2 [ 1/(1- u+ s/w(u)) - 1/(1- u)]) \to \exp(-s), \qquad u
\uparrow 1,$$ we have $F\in {\rm MDA}(\Lambda, w)$. Let $S\in (0,1)$
be a random variable independent of $R$ such that \eqref{eq:alp:den}
holds. Applying \netheo{T1} we obtain
\BQNY
\overline H(u) \sim \Gamma(\alpha+1)
\overline{F}(u)\overline{G}\left(1- \frac{(1- u)^2}{c_2u}\right), \quad u \uparrow 1.
\EQNY

\section{Applications}

\subsection{Ruin in the Presence of Risky Investments}

Consider the following discrete-time insurance risk model. Within
period $i$, the insurer's net profit (total premium income less
claim payment) is denoted by a real-valued random variable $Z_{i}$.
The insurer positions him/herself in a discrete-time financial
market consisting of a risk-free bond with a constant periodic
interest rate $\delta _{i}>0$ and a risky
stock with a periodic stochastic return rate $\Delta _{i}$ taking values in $%
(-1,\infty )$. Suppose that, in the beginning of each period $i$,
the insurer invests a fraction $\pi _{i}\in \lbrack 0,1)$ of his
current wealth in the stock and keeps the remaining wealth in the
bond. Denote by $U_{i}$
the insurer's wealth at time $i$, with a deterministic initial value $%
U_{0}=u\geq 0$. Then, $U_{i}$ evolves according to
\[
U_{i}=\left[ (1-\pi _{i})(1+\delta _{i})+\pi _{i}(1+\Delta
_{i})\right] U_{i-1}+Z_{i},\qquad i=1,2,\ldots .
\]

As usual, define the probability of ruin by time $n$ as
\[
\psi (u;n)=\pb{ \left. \min_{0\leq i\leq n}U_{i}<0\ \right\vert \
U_{0}=u} ,\qquad n=1,2,\ldots .
\]%
Assume that $Z_{1}$, $Z_{2}$, \ldots\ are independent and
identically distributed random variables, that $\Delta _{1}$,
$\Delta _{2}$, \ldots\ are
independent random variables, and that the two sequences $%
\{Z_{1},Z_{2},\ldots \}$ and $\{\Delta _{1},\Delta _{2},\ldots \}$
are mutually independent. Introduce
\begin{equation}
\Upsilon_i=\frac{1}{1+\Delta _{i}},\qquad R_{i}=-Z_{i},\qquad
S_{i}=\frac{1}{(1-\pi _{i})(1+\delta _{i})+\pi_i \left( 1+\Delta
_{i}\right) },\qquad i=1,2,\ldots . \label{two risks}
\end{equation}%
The random variable $\Upsilon_i$ is the random discount factor during period $%
i$ of the risky asset and it takes values in $(0,\infty )$, the
random
variable $R_{i}$ is the net loss during period $i$, and the random variable $%
S_{i}$ is the overall random discount factor during period $i$ of
the
investment portfolio. Denote by $F$ the common distribution of $%
\{R_{i},i=1,2,\ldots \}$.

According to Tang and Vernic (2010), if $F\in \mathcal{S}$, then
\begin{equation}
\psi (u;n)\sim \sum_{k=1}^{n}\pb{ R_{k}\prod_{i=1}^{k}S_{i}>u}
,\qquad u\rightarrow \infty . \label{estimate1}
\end{equation}%
See related discussions in Tang and Tsitsiashvili (2003, 2004). The
applicability of this formula requires explicit asymptotic
expressions for the tail probabilities in (\ref{estimate1}) and our
main results clearly are crucial for this purpose.

Notice from (\ref{two risks}) that if $\pk{ \Upsilon_i >u} $ is
regularly varying at infinity with index $-\alpha _{i}$ for some
$\alpha _{i}\in \lbrack 0,\infty ),$ then $\pk{
S_{i}>\hat{s}_{i}-1/u} $ is
regularly varying at infinity with index $-\alpha _{i}$, where $\hat{s}%
_{i}=(1-\pi _{i})^{-1}(1+\delta _{i})^{-1}\in (0,\infty )$.
Actually,
\begin{equation}
\pb{ S_{i}>\hat{s}_{i}-1/u} \sim \pb{ \Upsilon_i >\pi_{i} \hat{s}%
_{i}^{2}u} ,\qquad u\rightarrow \infty .  \label{iff}
\end{equation}%
Hence, if $F\in \mathcal{S}\cap \mathrm{MDA}(\Lambda ,w)$ and for each $%
1\leq i\leq n$, the tail probability $\pk{\Upsilon_i>u} $ is regularly
varying at infinity with index $-\alpha _{i}$ for some $\alpha
_{i}\in
\lbrack 0,\infty )$, then applying Theorem 3.1(a) to relation (\ref%
{estimate1}) we obtain, with $u_{k}=u\prod_{i=1}^{k}1/\hat{s}_{i}$,
\begin{eqnarray*}
\pb{ R_{k}\prod_{i=1}^{k}S_{i}>u} &=& \pb{
R_{k}\prod_{i=1}^{k}\frac{S_{i}}{\hat{s}_{i}}>u_{k}} \\
&\sim &\overline{F}\left( u_{k}\right) \prod_{i=1}^{k}\left[ \Gamma
\left(
\alpha _{i}+1\right) \pb{ \frac{S_{i}}{\hat{s}_{i}}>1-\frac{1}{%
u_{k}w(u_{k})}} \right]  \\
&\sim &\overline{F}\left( u_{k}\right) \prod_{i=1}^{k}\left[
\frac{\Gamma \left( \alpha _{i}+1\right) }{\left( \pi_{i}
\hat{s}_{i}\right) ^{\alpha _{i}}}\pk{ \Upsilon_i>u_{k}w(u_{k})}
\right] ,\qquad u\rightarrow \infty ,
\end{eqnarray*}%
where the last step is due to (\ref{iff}). We summarize all this as follows. \\

\BT \label{ruin}
Consider the discrete-time risk model introduced above. If $F\in \mathcal{S}%
\cap \mathrm{MDA}(\Lambda ,w)$ and for each $1\leq i\leq n$, the
tail probability $\pk{\Upsilon_i >u}$ is regularly varying at
infinity with index $-\alpha _{i}$ for some $\alpha _{i}\in \lbrack
0,\infty )$, then
\[
\psi (u;n)\sim \sum_{k=1}^{n}\overline{F}\left( u_{k}\right) \prod_{i=1}^{k}%
\left[ \frac{\Gamma \left( \alpha _{i}+1\right) }{\left( \pi_{i} \hat{s}%
_{i}\right) ^{\alpha _{i}}}\pk{ \Upsilon_i >u_{k}w(u_{k})} \right]
,\qquad u\rightarrow \infty .
\]
\ET Note in passing that if $F\in {\rm MDA}(\Lambda,w)$ with $w$
such that $\limit{u} w(u)=0$, then in order to show that $F\in
\kal{S}$ we can utilise \eqref{Tony}.

\subsection{Asymptotics of Conditional Tail Expectation}
Let $S_i$, $i\le n$, $R$ (with distribution function $F$), and $X_n$
be as in \netheo{T1}. If $F\in {\rm MDA}(\Lambda,w)$, then $w$
satisfies the self-neglecting property (see e.g., Reiss (1989) or
Resnick (2008))
 \BQN\label{eq:uv:b}
 \lim_{u \toomgF} \frac{ w(u+z/w(u))}{w(u)}=1
 \EQN
uniformly with respect to $z$ in every compact set of $\R$.  So if
the conditions of \netheo{T1}(a) are satisfied, it follows from
\eqref{eq:alp} and \eqref{eq:lem1:2}  that $H_n\in {\rm
MDA}(\Lambda,w)$. Consequently, we obtain the following asymptotic
formula (recall \eqref{eq:reco1})
 \BQN\label{eq:neu:0} \lim_{u
\uparrow r_F} \frac{ \E{ X_n- u\lvert X_n> u}}{\E{ R- u\lvert R>
u}}=1.
 \EQN
In several insurance and finance applications the mean
excess function is a crucial quantity (see Embrechts et al.\ (1997),
p. 294). The result in \eqref{eq:neu:0} shows that under the assumed conditions,    the
mean excess function is asymptotically invariant under random
contractions.

The Conditional Tail Expectation (CTE) for $R$ with continuous
distribution function is defined e.g., by
 $${\rm CTE}_R(u):= \E{R|R>u}= \E{ R- u\lvert R> u}+u, \qquad u>0.$$
Then in view of \eqref{eq:neu:0},
$$ \lim_{u \toomgF} \frac{ {\rm CTE}_{R}(u)}{u}=  \lim_{u \toomgF} \frac{ {\rm CTE}_{X_n}(u)}{u}=  1+\lim_{u \toomgF} \frac{ 1}{uw(u)}=1.$$
Consequently, if ${\rm VaR}_p(X_n)$ denotes the Value at Risk (VaR)
corresponding to the level  $p\in (0,1)$, i.e., $${\rm
VaR}_p(X_n):=\inf\left\{x:\pk{X_n \leq x}\geq p\right\},$$ then
$$ {\rm CTE}_{X_n}({\rm VaR}_{p}(X_n))\sim {\rm VaR}_{p}(X_n), \qquad p\uparrow 1.$$
It is well-known that for continuous risks CTE is more conservative
than VaR. The above asymptotics shows that in the Gumbel case CTE
and VaR are asymptotically the same and that this relation is
preserved under random scaling.

\subsection{Linear Combinations of Random Contractions}
In order to motivate the  next applications we consider a bivariate
scale mixture random vector $(U_1,U_2)$ with stochastic
representation \BQN\label{eq:spher} (U_1,U_2)\equaldis R (I_1 S,
I_2\sqrt{1- S^2}), \EQN
where $R$, with distribution function $F$, is almost
surely positive, $I_1,I_2$ assume values in $\{-1,1\}$, and $S$, with distribution function $G$, is a
scaling random variable taking values in $(0,1)$. Furthermore,
suppose that $I_1,I_2, R, S$ are mutually independent. If $S^2$
follows a beta distribution with parameters $1/2,1/2$ and
$\pk{I_1=1}=\pk{I_2=1}=1/2$, then $(U_1,U_2)$ is a spherically
distributed  random vector, see Cambanis et al.\ (1981). (We shorten
this by saying that $(U_1,U_2)$ is spherical).
 Hence by Lemma 6.1 of Berman (1983)
 \BQN\label{eq:berm61} c_1U_1+c_2U_2\equaldis \sqrt{c_1^2+c_2^2} U_1,
\qquad \forall c_1,c_2\inr.
 \EQN
If $S$ follows a beta distribution with parameters $a,b$, then
$(U_1,U_2)$ is a generalised Dirichlet random vector (see Hashorva
et al.\ (2007)), and \eqref{eq:berm61} does not hold in general.

Next, we derive the tail asymptotics of the aggregated risk
$$U(\rho)= \rho U_1+\sqrt{1-\rho^2}U_2, \qquad \rho\in (0,1)$$
for some general scaling random variable $S$. If $I_1=I_2=1$, then $U(\rho)$ is maximized with
 respect to $S$ at $S=\rho$. Hence we make the following assumption about the local form of $G$ at $\rho$:
\BQN\label{eq:G} \pk{\left\vert {S-\rho}\right\vert \leq
t}=L_\rho(t) t^{\alpha_\rho},\qquad \alpha_\rho\in [0,\IF),
 \EQN
for all $t\in (0,\ve)$, $\ve>0$, where $L_\rho$ is positive and
slowly varying at 0, and
 $L_\rho(0+)=0$ if  $\alpha_\rho=0$. Clearly, if $G$ possesses a density function $g$ continuous at $\rho$, then
\eqref{eq:G} holds  for any $\rho \in (0,1)$  with $\alpha_\rho=1$
and $L_\rho(t)= (2+o(1)) g(\rho)$ as $t\downarrow 0.$ \\

\begin{lem} \label{lem:11}
Let $I_1,I_2$ be two random variables taking values $-1,1$ with
$q=\pk{I_1=1,I_2=1}\in (0,1]$ and independent of a scaling random
variable $S$ with distribution function $G$. For given $\rho\in
(0,1)$ define a new random variable
$$S(\rho):= \rho I_1 S+ \sqrt{1- \rho^2} I_2  \sqrt{1-
S^2}.$$ If $G$ satisfies \eqref{eq:G} with $L_\rho(0+)=0$ when
$\alpha_\rho=0$, then
 \BQN \pk{S(\rho)> 1-
u}\sim q L_\rho(\sqrt{u}) (2 u  (1- \rho ^2))^{\alpha_\rho/2},
\qquad u\downarrow 0.
 \EQN
\end{lem}

Note that if $G$ is absolutely continuous with a positive density
function $g$ continuous at $\rho$, then we have
 \BQN \pk{S(\rho)> 1-
u}\sim q g(\rho) \sqrt{8  (1- \rho ^2) u} , \qquad u\downarrow 0.
 \EQN

In view of \netheo{T1} and \nelem{lem:11}, we now study the tail
behavior of $U(\rho)$. First consider the \emph{Gumbel} case i.e.,
$F\in {\rm MDA}(\Lambda,w)$ with $r_F\in (0,\IF]$. Then with
$\eta(u)= u w(u)$ for $u>0$ we have
 \BQN\label{eq:UR:G} \pk{U(\rho)> u}\sim
q\Gamma(\alpha_\rho/2+1) L_\rho(\eta(u)^{-1/2}) \fracl{2  (1- \rho
^2)}{\eta(u)}^{\alpha_\rho/2}
 \overline F(u), \qquad u\uparrow r_F.
 \EQN
Since $\lim_{u\uparrow r_F} \eta(u)= \IF$ it follows that
$\lim_{u\uparrow r_F}\pk{U(\rho)>u}/\overline F(u)=0.$ If $G$
possesses a density function $g$ continuous at $\rho$, then
\BQN\label{eq:UR:W} \pk{U(\rho)> u}\sim q \Gamma(1/2) g(\rho)
\fracl{2  (1- \rho ^2)}{\eta(u)}^{1/2}
 \overline F(u), \qquad u\uparrow r_F.
\EQN

Now if $(U_1,U_2)$ is spherical, then \eqref{eq:berm61} (obviously!)
implies the tail equivalence
 \BQN\label{eq:condF1} \pk{U(\rho)>
u}\sim \pk{U_1>u}, \qquad u\uparrow r_F.
 \EQN
The corresponding density function $h_\rho$  of $U(\rho)$ satisfies
 \BQN\label{eq:condF1:de} h_\rho (u)\sim h_1(u), \qquad u\uparrow
r_F.
 \EQN

The following converse result characterises spherical random
vectors. \\

\BT \label{theo:chr1} Suppose that $I_1,I_2$ are independent random
variables assuming values $-1,1$ with probability $1/2$, the
distribution function $F$ with $F(0-)=0$ satisfies \eqref{eq:LL},
and $G$ possesses a continuous density function $g$.  Then
\eqref{eq:condF1} holds for all $\rho \in (0,1)$ if and only if
$(U_1,U_2)$ is  spherical. Similarly,  \eqref{eq:condF1:de} holds
for all $\rho \in (0,1)$ if and only if $(U_1,U_2)$ is spherical.
\ET

Next consider the \emph{Weibull} case. Assume that $F\in {\rm
MDA}(\Psi_\gamma)$ with $\gamma \in (0,\IF)$ (and $r_F=1$). Applying
\netheo{T1} and \nelem{lem:11} we obtain
 \BQNY \pk{U(\rho)> 1- u}\sim q
\frac{\Gamma(\alpha_\rho/2+1) \Gamma(\gamma+1)}
{\Gamma(\alpha_\rho/2+\gamma+1)} L_\rho(u^{1/2}) (2  u (1- \rho
^2))^{\alpha_\rho/2} \overline F(1-u), \qquad u\downarrow 0.
 \EQNY
If $G$ possesses a continuous density function $g$, then this
simplifies to
 \BQNY \pk{U(\rho)> 1- u}\sim q \frac{\Gamma(1/2)
\Gamma(\gamma+1)} {\Gamma(3/2+\gamma)}  g(\rho) (2 u (1- \rho
^2))^{1/2}  \overline F(1-u), \qquad u\downarrow 0.
 \EQNY

The following result gives the Weibull analogue of
\netheo{theo:chr1}. \\

\BT \label{chr2} Suppose that $I_1,I_2$, $g,G$ are as in
\netheo{theo:chr1}.  If $F\in {\rm MDA}(\Psi_\gamma)$ with $\gamma
\in (0,\IF)$ (and $r_F=1$), then \eqref{eq:condF1} holds for all
$\rho \in (0,1)$ if and only if $(U_1,U_2)$ is  spherical.
 \ET
We remark that if $F\in {\rm MDA}(\Phi_\gamma)$ with $\gamma \in
(0,\IF)$, then the tail behaviour of $U(\rho)$ follows from
Breiman's lemma. Indeed, under this assumption $U(\rho)$ has a
distribution function in ${\rm MDA}(\Phi_\gamma)$.

\subsection{Max-Domain of Attraction of Bivariate Samples}
Suppose that $Q$ is the distribution function of a bivariate random
vector $(X,Y)$. Extending \eqref{eq:LL}, we say that $Q$ belongs to
the max-domain of attraction of a bivariate max-stable distribution
function $N$ if
 \BQN \label{eq:LLL}
 \limit{n}\sup_{x,y\inr} \Abs{Q^n(a_nx+ b_n,c_ny+ d_n)- N(x,y)}= 0
 \EQN
holds for some constants $a_n>0,c_n>0, b_n,d_n\inr,n\ge 1$. This
implies that each univariate marginal distribution of $Q$ is in the
max-domain of attraction of the corresponding univariate marginal of
$N$. Conversely, if further
$$ \limit{n} n\pk{X>b_n, Y>d_n}=0,$$
then  \eqref{eq:LLL} holds with $N$ a product distribution function
with univariate extreme value marginal distributions.

In this last application we discuss the  max-domain of attraction
for the distribution function of $(U_1, U(\rho))$, defined in the
previous subsection, assuming that $F$ is in the max-domain of
attraction of some univariate distribution. For any $\rho\in
(0,1)$ denote by $Q_\rho$ the distribution function of the bivariate
random vector $(U_1, U(\rho))$, and let $Q_{i,\rho}$, $i=1,2,$
denote the corresponding marginal distribution functions. We focus
here on the cases where $F$ is in either the Gumbel or the Weibull
max-domain of attraction. Assume that the distribution function $G$
of $S$ satisfies \eqref{eq:G} with $\alpha_\rho\in [0,\IF)$, and
$L_\rho$ positive, slowly varying at 0, and $L_\rho(0+)=0$ if
$\alpha_\rho=0$.

First consider the \emph{Gumbel} case.  Assume that $F\in {\rm
MDA}(\Lambda,w)$. It follows from \eqref{eq:UR:G}  that $Q_{i,\rho}
\in {\rm MDA}(\Lambda,w)$ for $i=1,2$. We show that, if $\rho\in
[0,1)$, then  $Q_\rho$ is in the max-domain of attraction of a
bivariate distribution function which is a product distribution.  If
$r_F$ is finite, then this follows from the evident fact that  $U_1$
and $U(\rho)$ cannot simultaneously take their maximum values. If
$r_F=\IF$, then set $b_i(n):= Q_{i,\rho}^{-1}(1- 1/n), n>1$ with
$Q_{i,\rho}^{-1}$ the generalized inverse of the $i$-th marginal
distribution of $Q_\rho$. It follows from Lemma 4.2 of Hashorva and
Pakes (2010) that
 \BQN\label{eq:con:b2b1} \limit{n} (b_2(n)- \rho b_1(n)) = \IF.
 \EQN
Our claim follows if \BQN \label{suff}
 \limit{n} n\pk{ U_1> b_1(n), U(\rho)> b_2(n)} = \limit{n} \pk{ U(\rho)> b_2(n)\lvert U_1> b_1(n)}=0.
 \EQN
In view of Theorem 2 of Hashorva (2009)
$$ \limit{n} \pb{ U(\rho)> \rho b_1(n)[1+ x/\sqrt{ b_1(n) w(b_1(n))}] \bigl \lvert U_1> b_1(n)}\in(0,\IF), \quad \forall
x\inr$$
implying \eqref{suff}, and hence our claim.
Note that condition \eqref{eq:con:b2b1} is satisfied for a
distribution function $F$ with tail asymptotics given by
\eqref{eq:kotz:Fk}.

Next consider the \emph{Weibull} case. Assume that $F\in {\rm
MDA}(\Psi_\gamma)$ with $\gamma\in (0,\IF)$ and $r_F=1$. In view of
\eqref{eq:UR:G} it follows that $Q_{1,\rho} \in {\rm
MDA}(\Psi_{\gamma+\alpha_1/2})$ and $Q_{2,\rho} \in {\rm
MDA}(\Psi_{\gamma+\alpha_\rho/2})$.  Since $r_F=1$ it follows that
$Q_\rho$ is in the max-domain of a bivariate distribution function
which is a product distribution. By \nelem{000}, this outcome can be formally generalized to the case where $r_F$ is
an arbitrary positive constant.

\section{Proofs}

\prooflem{000} In Lemma 5.2 of Hashorva and Pakes (2010) the proof
of \nelem{000} is shown for $F\in {\rm MDA}(\Lambda,w)$ and $c=1$.
The case $c\in (0,\IF)$ and $F$ is in the \FRE or Weibull max-domain
of attraction can be easily shown, therefore omitted here. \QED

\prooflem{interesting} First of all, if relation \eqref{Tony} holds
for some $\lambda >0$, then it
holds for all $\lambda >0$. This can be verified as follows. For all $%
0<\lambda _{1}<\lambda _{2}<\infty $, by $F\in \mathrm{MDA}(\Lambda
,w)$ we have, as $u\rightarrow \infty $,
\begin{eqnarray*}
\int_{u-\lambda _{2}/w(u)}^{u-\lambda
_{1}/w(u)}\overline{F}(u-y)dF(y) &\leq
&\overline{F}(\lambda _{1}/w(u))\left[ \overline{F}(u-\lambda _{2}/w(u))-%
\overline{F}(u-\lambda _{1}/w(u)))\right] \\
&=&o\left( (e^{\lambda _{2}}-e^{\lambda _{1}})\overline{F}(u)\right)
,
\end{eqnarray*}%
and similarly, $\int_{\lambda _{1}/w(u)}^{\lambda _{2}/w(u)}\overline{F}%
(u-y)dF(y)=o\left( \overline{F}(u)\right) $.

Next, we assume $F\in \mathcal{L}\cap \mathcal{S}$. Since $%
\lim_{u\rightarrow \infty }w(u)=0$, it holds for every $\lambda >0$, every $%
T>0$, and all large $u>0$ that
$$
\int_{\lambda /w(u)}^{u-\lambda /w(u)}\overline{F}(u-y)dF(y)\leq \overline{%
F^{\ast 2}}(u)-\int_{-\infty }^{T}\overline{F}(u-y)dF(y)-\int_{u+T}^{\infty }%
\overline{F}(u-y)dF(y).
$$%
By $F\in \mathcal{L}\cap \mathcal{S}$, it holds that $\overline{F^{\ast 2}}%
(u)\sim 2\overline{F}(u)$, that \BQN \label{a} \int_{-\infty
}^{T}\overline{F}(u-y)dF(y)\sim \overline{F}(u)F(T),
 \EQN
and that
 \BQN \label{b}
\overline{F}(-T)\leq \liminf_{u\rightarrow \infty }{\frac{1}{\overline{F}(u)}%
}\int_{u+T}^{\infty }\overline{F}(u-y)dF(y)\leq
\limsup_{u\rightarrow \infty
}{\frac{1}{\overline{F}(u)}}\int_{u+T}^{\infty
}\overline{F}(u-y)dF(y)\leq 1. \EQN Then, relation \eqref{Tony}
follows since $T$ can be arbitrarily large.

Finally, we assume that relation \eqref{Tony} holds for all $\lambda
>0$ and we prove $F\in \kal{S}$. This part of Lemma \ref{interesting} is an easy consequence of
Theorem 3.6 of  Foss et al.\ (2009). Actually, for a distribution
$F\in {\rm MDA}(\Lambda,w)$ with $r_F=\IF$ and $\limit{u} w(u)=0$,
for any $\lambda_u>0$ with $\limit{u} \lambda_u =0$ the function
$h(u)=\lambda_u/w(u)$ is $F$-insensitive in the sense of Foss et
al.\ (2009), i.e.,
$$ \limit{u} \frac{ \overline{F}(u+ h(u))}{\overline{F}(u)}=1.$$
Hence by Theorem 3.6 of  Foss et al.\ (2009), relation \eqref{Tony}
implies $F\in \kal{S}$. Nevertheless, we give here another
self-contained proof. Similarly as above, for every $\lambda
>0$ and $T>0$ we write
$$
\overline{F^{\ast 2}}(u)=\left( \int_{-\infty
}^{T}+\int_{T}^{\lambda /w(u)}+\int_{\lambda /w(u)}^{u-\lambda
/w(u)}+\int_{u-\lambda /w(u)}^{u+T}+\int_{u+T}^{\infty }\right)
\overline{F}(u-y)dF(y).
$$%
Estimates of the first and last terms above have been given in
\eqref{a} and \eqref{b}, and an estimate of the third term is given by
(\ref{Tony}). For the second and the fourth terms, we have, as
$u\rightarrow \infty $,
$$
{\frac{1}{\overline{F}(u)}}\int_{T}^{\lambda /w(u)}\overline{F}%
(u-y)dF(y)\leq {\frac{\overline{F}(u-\lambda
/w(u))}{\overline{F}(u)}}\left(
\overline{F}(T)-\overline{F}(\lambda /w(u))\right) \rightarrow e^{\lambda }%
\overline{F}(T),
$$%
and
$$
{\frac{1}{\overline{F}(u)}}\int_{u-\lambda /w(u)}^{u+T}\overline{F}%
(u-y)dF(y)\leq {\frac{\overline{F}(u-\lambda /w(u))-\overline{F}(u+T)}{%
\overline{F}(u)}}\rightarrow e^{\lambda }-1.
$$%
By the arbitrariness of $\lambda $ and $T$, we easily conclude that $%
\overline{F^{\ast 2}}(u)\sim 2\overline{F}(u)$.\QED

\prooftheo{T1}  a) The claim in \eqref{eq:lem1:1} follows by Lemma
A.5 of Tang and Tsitsiashvili (2004) and Theorem 3 of Hashorva
(2009). We give here another direct proof. Let $X=RS$, where the
factors are independent and the distribution function $G$ of  $S$
satisfies \eqref{eq:alp} with the subscript $i$ omitted, i.e., \BQN
\label{RVG} \overline G(1-x)=x^\alpha L(1/x), \EQN
with  $\alpha\ge
0$ and $L$ slowly varying at infinity.  Clearly,
\BQN\label{eq:di:H} \pk{X > u}= \int_u^ {r_F} \overline G(u/y)\,
dF(y).
\EQN
Substitute $y=u+z/w(u)$ and define the random variable
$W_u$ by
$$ \pk{W_u>z}={\overline F(u+z/w(u)) \over \overline F(u)},\qquad 0\le z<r(u),$$
where $r(u)=(r_F-u)w(u)<\infty$ if $r_F< \IF$, and $r(u)=\IF$
otherwise. Observe that \eqref{eq:rdfd} can be expressed as $W_u
\st{d}\to W, u \uparrow r_F$, where $W$ is a random variable
following an exponential distribution with mean $1$ and $\st{d} \to$ means convergence in distribution.

Next, \eqref{eq:di:H} can be expressed as
$$\pk{X>u}=\overline F(u) \E{\overline{G}(1/(1+W_u/\eta(u)))}, \qquad \eta(u)=uw(u).$$
The survival function in this expectation is asymptotically
proportional to $W_u^\alpha \overline G(1-1/\eta(u))$ almost surely
as $u \uparrow r_F$.
We consider two cases.\\
i) If $W_u\le 1$, then $W_u/(W_u+\eta(u))\le 1/(1+\eta(u))$, and hence dominated convergence implies that
$$\lim_{u\uparrow r_F}{ \E{\overline G(1-W_u/(W_u+\eta(u))); W_u\le 1} \over \overline G(1-1/(1+\eta(u)))}=\E{W^\alpha; W\le
1}.$$
ii) If $W_u\ge 1$, then Potter's bounds for slowly varying functions
(Bingham et al. (1987, Theorem 1.5.6 (i))) imply that for chosen
constants $A>1$ and $\delta>0$, there exists a number $r>0$ such
that, once $\eta(u)>r$ (see \eqref{eq:uv}),
$${L(1+\eta(u)/W_u) \over L(1+\eta(u))} \le AW_u^\delta.$$
Let  $M_n$ denote the maximum of $n$ independent copies of $R$. It follows from \eqref{eq:LL} that, for all positive $k$,  the
$k$-th order moments of $|M_n-b_n|/a_n$ converge to the $k$-th order moment of the Gumbel distribution (Pickands (1968)). This
can be utilised to show that
  $\{W_u^{\alpha+\delta}; u<r_F\}$ is a uniformly integrable family. So it follows from \eqref{RVG} that
$$\lim_{u \uparrow r_F} { \E{\overline G(1-W_u/(W_u+\eta(u))); W_u\ge 1} \over \overline G(1-1/(1+\eta(u)))}=\E{W^\alpha; W\ge
1}.$$
Since $\E{W^\alpha}=\Gamma(1+\alpha)$, we conclude that, as
$u\uparrow r_F$,
$$\pk{X>u}\sim \Gamma(1+\alpha)\overline G\left(1-1/\eta(u)\right) \overline F(u).$$
The self-neglecting property of the scaling function (see
\eqref{eq:uv:b}) implies in addition that $H\in {\rm
MDA}(\Lambda,w)$. Hence \eqref{eq:lem1:1} can be proved by
induction.

a1) The proof follows directly from \nelem{interesting}.

b) Since $r_F=1$, we express \eqref{eq:di:H} as
$$\pk{X>1-u}=\int_{1-u}^1 \overline G((1-u)/y)\, d F(y)=\overline F(1-u)\eg{G\Bigl( \frac{1-u}{1-uW_u}\Bigr)},$$
where
$$\pk{W_u\le z}={\overline F(1-uz) \over \overline F(1-u)},\qquad 0\le z\le 1.$$
So as $u\to0$, we have $W_u\st{d}\to W$, where $\pk{W\le
z}=z^\gamma$ for $0\le z\le 1$. Since
$${1-u \over 1-uW_u}=1-{u(1-W_u) \over 1-uW_u},$$
it follows from \eqref{RVG} that, almost surely,
$$\overline G\left({1-u \over 1-uW_u}\right)\sim (1-W_u)^\alpha \overline G(1-u), \qquad u\to0.$$
Clearly, for any $u\in (0,1)$
$$ \overline G\left({1-u \over 1-uW_u}\right) \le \overline G(1-u),$$
hence dominated convergence yields that
$$\lim_{ u\downarrow 0} \frac{1}{\overline G(1-u) } \eg{G\Bigl( \frac{1-u}{1-uW_u}\Bigr)} = \E{(1-W)^\alpha}={\Gamma(\alpha+1)\Gamma(\gamma+1) \over \Gamma(\alpha+\gamma+1)}.$$

We prove  \eqref{eq:lem1:2} by induction as follows.  Let
$$C_n={\Gamma(\gamma+1) \over \Gamma(\gamma_n+1)}\prod_{i=1}^n \Gamma(\alpha_i+1), \qquad \gamma_n= \gamma+
\sum_{i=1}^n \alpha_i,$$
and assume that
$$\overline H_{n-1}(1-u)\sim C_{n-1} \left[\prod_{i=1}^{n-1}\overline G_i(1-u)\right] \overline F(1-u).$$
Noting that the right-hand side is regularly varying at zero with
index $\gamma_{n-1}$, the case $n=1$ implies that
$$\overline H_n(1-u)\sim C_{n-1}\eg{(1-W_{(n-1)})^{\alpha_n}}
 \left[\prod_{i=1}^{n} \overline G_i(1-u) \right] \overline F(1-u),$$
 where $\pk{W_{(n-1)}\le z}=z^{\gamma_{n-1}}$. But
 $$C_{n-1}\eg{(1-W_{(n-1)})^{\alpha_n}}=C_{n-1}
 {\Gamma(\alpha_n+1)\Gamma(\gamma_{n-1}+1) \over \Gamma(\alpha_{n}+\gamma_{n-1}+1)}=C_n,$$
 and the assertion follows.\QED


\prooftheo{T2} In all cases it suffices to prove the case $n=1$
because the general result follows from asymptotic estimates
obtained from \netheo{T1}  applied to the product with $n-1$
contraction factors. Thus, omitting subscripts, the density function
of $X=RS$ is (see e.g., Lemma 2.1 of Pakes and Navarro (2007))
 \BQN
\label{pdf} h(u)=\int_u^{r_F} y^{-1}g(u/y)dF(y)= \int_0^1 y^{-1}
f(u/y)dG(y),
 \EQN
 where the first equality applies when $G$ has a
density function, and the second equality applies when $F$ has a
density
function.\\

a1) Let $W_u$ be a random variable whose distribution function is
$$\pk{W_u\le z}={\overline F(u/z) \over \overline F(u)},\qquad z \in (0,1].$$
Clearly $W_u \st{d}\to W$ as $u \to \infty$, where $\pk{W\le
z}=z^\gamma$ (so it is degenerate at $0$ provided that $\gamma=0$).
Substituting $y=u/z$ in the first integral of \eqref{pdf} yields
$$h(u)=\frac{\overline F(u)}{u} \E{W_ug(W_u)}.$$
It follows from dominated convergence that the expectation converges
to
$$\lim_{u\to\infty}\E{W_ug(W_u)}=\E{Wg(W)}=\gamma \int_0^1 g(z)z^\gamma dz=\gamma \E{S^\gamma}, $$
and hence, from \eqref{eq:breiman}, that
 \BQN \label{bdf}
\lim_{u\to\infty} {u h(u) \over \overline H(u)}=\gamma.
 \EQN

 a2) The second integral form in \eqref{pdf} can be expressed as $h(u)=\E{S^{-1}f(uS^{-1})}$. Expressing the regular
 variation assumption  as $f(u)=u^{-\gamma-1}L(u)$, where $L$ is slowly varying at infinity, we have
 $${h(u) \over f(u)}=\E{S^\gamma [L(uS^{-1})/L(u)]}.$$
Choose $\epsilon \in(0,\gamma)$ if $\gamma>0$ or choose as in the
assumptions if $\gamma=0$. Since $S\le1$,
 it follows from Potter's bounds that for chosen $A>1$ there exists $u'>0$ such that  $L(uS^{-1})/L(u) \le S^{-\epsilon}$
 for all $u>u'$.
 Hence dominated convergence yields $\limit{u}h(u)/f(u)= \E{S^\gamma}$.
  The Karamata-Abelian theorem (Bingham et al. (1987, p. 26)) implies that \eqref{bdf} still holds.

b) As in the proof of \netheo{T1} it suffices to consider the case
$n=1$, i.e., $X=RS$, where the density function $g$ of $S$ has the
form
 \BQN \label{grv}
 g(1-1/u)=\alpha u^{-\alpha+1}L(u), \qquad u>1.
 \EQN

Letting $\eta(u)=uw(u)$, and recalling notation from the proof of
\netheo{T1}(a), it follows from \eqref{pdf} that the density
function of $X$ satisfies
$$h(u)\sim u^{-1} g(1-1/\eta(u))\overline F(u) \E{(1+W_u/\eta(u))^{-1}R_u}, \qquad u \uparrow r_F,$$
where
\begin{eqnarray*}
R_u &=&  {g(1-W_u/(\eta(u)+W_u)) \over g(1-1/(\eta(u)+1))}= W_u^{\a-1} \left({1+\eta(u) \over W_u+\eta(u)}\right)^{\a-1} {L(1+\eta(u)/W_u) \over L(1+\eta(u))}.
\end{eqnarray*}
Note that $(1+W_u/\eta(u))^{-1}<1$ and that the second and third
factors of $R_u$ converge to unity almost surely.

Let $l$ be a (large) positive constant and  $B_u=\{W_u\le \eta(u)/l\}$. Note that $\pk{B_u}\to1$ as $u\uparrow r_F$. If
$\alpha\ge 1$, then the second factor of $R_u$ is bounded above by
$(1+1/\eta(u))^{\alpha-1}\to1$, and if $0<\alpha<1$, then on $B_u$  the second factor equals
$$\left({W_u+\eta(u) \over 1+\eta(u)}\right)^{1-\alpha}\le
\left({(1+l^{-1})\eta(u) \over 1+\eta(u)}\right)^{1-\alpha}
\le\left(1+l^{-1}\right)^{1-\alpha}.$$

Next, choosing $\delta\in(0,\max(\alpha,1))$ and $A>1$, $l$ can be
made so large that on $B_u$ and with $u$ such that $\eta(u)>l$, the
third factor of $R_u$ is dominated by Potter's bound
$A\max(W_u^\delta,W_u^{-\delta})$.  Dominated convergence hence
gives the conclusion
$$\lim_{u\toomgF} \E{(1+W_u/\eta(u))^{-1}R_u; W_u\le \eta(u)/l}=\E{W^{\alpha-1}}.$$

Assume that $\eta(u)\ge l$ and consider outcomes on $\overline B_u$,
i.e., that $l \le \eta(u) < lW_u$. If $\alpha\ge 1$, then the second
factor of $R_u$ is bounded above by unity, and if $0<\alpha<1$ then
this factor is dominated by $W_u^{1-\alpha}$. To deal with the third
factor, observe that condition \eqref{gbb} is equivalent to the
existence of $p>0$ such that $L(1+z)=O(z^{-p})$ as $z\to0$. Consequently,
$l$ can be chosen so large that the numerator in the third factor of
$R_u$ is $O_p[(W_u/\eta(u))^p]$, as $u\toomgF$. It follows from
these estimates that
$$\E{R_u; \overline B_u} = O\left( { \E{W_u^{p+\max(\alpha-1,0)};\overline B_u} \over (\eta(u))^pL(\eta(u))}\right)\to 0, \qquad u\toomgF$$
since the denominator tends to infinity.

Since $ \E{W^{\alpha-1}}=\Gamma(\alpha)$, it follows that,
$$h(u)\sim \Gamma(\alpha) u^{-1}g(1-1/\eta(u))\overline F(u) .$$
But $g(1-y)\sim (\alpha/y)\overline G(1-y)$, so \eqref{eq:T2:1}
(with the subscripts omitted) follows. The general result follows
from \netheo{T1} after replacing $F$ with the distribution function
of $R\prod_{j\ne i}S_j$.

c) The proof for the case $n=1$ is similar to those above. Observe
that  if $u,w\in(0,1)$, then $(1-uw)/(1-w)>1$. Assuming \eqref{grv},
then choosing $\delta\in(0,\alpha)$, we have Potter's bound
$${L((1-uW_u)/(u(1-W_u))) \over L(1/u)}\le A\max(1, (1-uW_u)^\delta(1-W_u)^{-\delta}),\qquad u>u',$$
where $W_u$ is as in the proof of \netheo{T1}(b). This leads to the asymptotic form
$$h(1-u)\sim g(1-u)\overline F(1-u) \E{(1-W)^{\alpha-1}}, \qquad u \downarrow 0.$$
The assertion for $n=1$ then follows with no further assumptions about $g$.

For the general case we may assume with no loss of generality that \eqref{grv} holds with $i=n$, and then apply the single
factor case with $H_{n-1}$ replacing $F$. It follows from \netheo{T1}(b) that
$$h_n(1-u)\sim g_n(1-u)\overline H_{n-1}(1-u)\E{(1-W_{n-1})^{\alpha_n-1}}.$$
But $g_n(1-u)\sim (\alpha_n/u)\overline G_n(1-u)$ and
$$\alpha_nC_{n-1}\E{(1-W_{n-1})^{\alpha_n-1}}=C_{n-1}
{\alpha_n \Gamma(\alpha_n) \Gamma(\gamma_{n-1}+1) \over
\Gamma(\gamma_n)}=\gamma_nC_n.$$ It follows that $h_n(1-u)\sim
\gamma_n \overline H_n(1-u)/u$, whence  the general result
\eqref{eq:T2:2}. \QED

\prooftheo{TT} It suffices to consider the case $n=1$. We can, with no loss in generality, choose the scaling function to be
differentiable and satisfy
$$\overline w(s) := {d \over ds} {1\over w(s)}\to 0, \qquad u\uparrow r_F.$$
Replacing $u$ in \eqref{pdf} with $u+x/w(u)$, the substitution $y=s+x/w(s)$ yields
$$h(u+x/w(u))=\int_u^{r_F} g\left({u+x/w(u) \over s+x/w(s)}\right)
{1+\overline w(s) \over s+x/w(s)} f(s+x/w(s)) ds.$$ Write the
argument of $g$ as $(u/s)a(u,s)$, where
$$a(u,s)={1+x/(u w(u)) \over 1+x/(s w(s))}\to 1,\qquad  u\uparrow r_F.$$
Also, the middle factor in the above integrand is asymptotically
proportional to $s^{-1}$. Hence \eqref{eq:vm:h} follows from
\eqref{eq:vm:f} and \eqref{eq:gu}; see Lemma 5.1 of  Takahashi and
Sibuya (1998).\QED

\prooflem{lem:11} Some algebra shows that $S(\rho)\le 1$ and it is
bounded away from unity unless $I_1=I_2=1$. If this event occurs,
then $S(\rho)$ is close to $1$ if and only if $S$ is close to
$\rho$. So on the event $\{I_1=I_2=1\}$, algebra reveals that, for
all small $u>0$,  $ S(\rho)>1-u$ if and only if
$$(\rho-S)^2+2\rho Su<2u-u^2,$$
and this condition is equivalent to
$$(\rho-S)^2<(1+o_p(1))2(1-\rho^2)u.$$
But
$$\pk{(\rho-S)^2\le t} \sim t^{a_\rho/2}L_{\rho}(\sqrt t), \quad u \downarrow 0$$
and the assertion follows. \QED

\prooftheo{theo:chr1} The assumptions imply that \eqref{eq:UR:W}
holds for all $\rho\in(0,1)$. The tail equivalence \eqref{eq:condF1}
implies that the factor $g(\rho)(1-\rho^2)^{1/2}$ is constant with
respect to $\rho\in(0,1)$. The density function of $S^2$ is
$g(y^{1/2})/2y^{1/2}$, i.e., $S^2$ has the beta distribution with
parameters $(1/2,1/2)$. If \eqref{eq:condF1:de} holds, then so does
\eqref{eq:condF1}.  Hence either condition implies that $(U_1,U_2)$
is spherical.  \QED

\prooftheo{chr2} The proof is similar to  that of \netheo{theo:chr1}. \QED


\baselineskip14pt

\bibliographystyle{plain}

\end{document}